\documentclass[a4paper,notitlepage, 8pt,reqno]{article}
  \usepackage{graphicx}
  \usepackage{mathtext}
 \usepackage[cp1251]{inputenc}
  \usepackage[T2A]{fontenc}
 \usepackage[russian,english]{babel}

 \usepackage[all,arc,poly, 2cell,curve,arrow,tips]{xy}

  \usepackage{enumerate, eucal, amsthm,amsmath, amssymb}

  \allowdisplaybreaks[4]

 \theoremstyle{definition}
 \newtheorem{defn}{Definitio}%[section]

 \theoremstyle{plain}
 \newtheorem{thm}{Theorem}
 \newtheorem*{thm*}{Theorem}
 \newtheorem{prop}{Proposition}
  \newtheorem*{prop*}{Предложение}
 \newtheorem{cor}{Corollary}
  \newtheorem*{cor*}{Следствие}
 \newtheorem{lem}{Lemma}
  \newtheorem*{lem*}{Лемма}

 \theoremstyle{remark}

 \newtheorem*{remark*}{Замечание}

 \renewcommand{\abstractname}{}

  \newcounter{ab}
\title{The comparison Gelfand-Tsetlin-Molev and  Gelfand-Tsetlin-Zhelobenko bases for $\mathfrak{sp}_{2n}$ }

 \author{D. V. Artamonov\footnote{Lomonosov Moscow State University, artamonov.dmitri@gmail.com}}
\begin{document}
 \maketitle

\renewcommand{\abstractname}{}

\begin{abstract}
А construction of Gelfand-Tsetlin type base vectors  in a  finite-dimensional representation of  $\mathfrak{sp}_{2n}$ was firstly obtained  in 60-th by Zhelobenko \cite{zh1}.  But the final construction was obtained only in the year 1998 by Molev \cite{M1}, who gave a construction of  Gelfand-Tsetlin type base vectors and derived formulas for the action of generators of the algebra in this base. These two approaches use different ideas.

 In the present paper we compare these two approaches. Also we show that the Molev's base vectors can be obtained using a construction based on a relation between restriction problems $\mathfrak{sp}_{2n}\downarrow\mathfrak{sp}_{2n-2}$  and $\mathfrak{gl}_{n+1}\downarrow\mathfrak{gl}_{n-1}$, analogous  to the construction giving the Zhelobenko's base.

\end{abstract}

%\tableofcontents

\section{Introduction}

In the year
 1950 there appeared two short articles of Gelfand and Tsetlin
\cite{GC1}, \cite{GC2},  where they gave explicit construction of finite-dimensional representation for the Lie algebras $\mathfrak{gl}_n$ and
$\mathfrak{o}_n$. More precise they gave indexation of base vectors and obtained  formulas for the action of generators in the base. In these notes there were not given a derivation of the presented formulas. In the case of $\mathfrak{gl}_n$ the derivation was published by Zhelobenko in the paper \cite{zh1} (see also
\cite{zh2}).

A natural question is the construction of an analogue of the Gelfand-Tsetlin base for sympletic algebras. However the attempt to apply the ideas of  Gelfand-Tsetlin in this situation meets a new phenomenon. The construction a base for the series of algebras
$g_n$  (where $g_n=\mathfrak{gl}_{n}$ or $\mathfrak{sp}_{2n}$)
is based on investigations of the branching of an irreducible representation under the restriction of the algebra $g_n\downarrow g_{n-1}$. In the case
$\mathfrak{gl}_{n}$ the spectrum of possible
$g_{n-1}$-irreps is simple and in the case of
$\mathfrak{sp}_{2n}$ can have multiplicities.

To construct a Gelfand-Tstlin type base for the symplectic algebra one has to find a way of construction of a base in the space of
$\mathfrak{sp}_{2n-2}$-highest vectors in a given irrep of $\mathfrak{sp}_{2n}$.

For the first time the solution of this problem was given by Zhelobenko in \cite{zh1} (see also
\cite{zh2}).  However he managed not to find the formulas for the action of generators in his base.

The full solution of the problem of construction of the
Gelfand-Tsetlin type base for the algebra $\mathfrak{sp}_{2n}$
(including the derivation of the formula for the action of
generators of the algebra in this base)   was given by Molev in
\cite{M1} (see also \cite{M}).  He used a new technique: the main
step in the construction of the base  was a construction of the
action of the yangian  $Y(\mathfrak{gl}_2)$ on the space of
$\mathfrak{sp}_{2n-2}$-highest vectors with a fixed
$\mathfrak{sp}_{2n-2}$-weight.

Mention that the indexation of the base vectors in the approaches of Molev and Zhelobenko is the same.

Thus there naturally appears a question: do these constructions define the same base of not.
The answer to this question is negative. The reason is simple: these two construction use different subalgebras
$\mathfrak{sp}_{2n-2}\subset\mathfrak{sp}_{2n}$.  More precise, let $\mathfrak{sp}_{2n}$ be
a Lie algebra that preserves the symplectic form $\sum_{i=1,...,n}dx_{-i}\wedge dx_{i}$ in the space
 $\mathbb{C}^{2n}$ with coordinates  $x_{-n},...,x_{n}$. Then Zhelobenko takes the subalgebra  $\mathfrak{sp}_{2n-2}$
 that preserves the coordinates  $x_{-n},...,x_{-2},x_{2},...,x_n$\footnote{When we say that a subalgebra preserves
  some coordinates we mean that the subspace
  spanned by these coordinates is invariant and onto complement coordinates the algebra act as zero},
   and Molev takes a subalgebra that preserves the coordinates $x_{-n+1},...,x_{n-1}$. Thus the Gelfand-Tsetlin type bases are different.

But one can ask then the following question. What are the properties
of the automorphism of an $\mathfrak{sp}_{2n}$-irrep that sends a vector of the base of
 Zhelobenko to the vector of the base of Molev corresponding to the same diagram. The answer is
 given in the Theorem \ref{osnvteor}. It turns out that this mapping preserves modulo multiplication by a constant the action
 of the so called  raising and lowering operators  that form the Mikelsson-Zhelobenko algebra (see it's definition in Section \ref{defmz}).

As a byproduct we obtain isomorphisms between different
Mickelsson-Zhelobenko algebras (see Lemma \ref{soo}). Also we obtain an important result saying that the Molev's base, that was originally obtained using a construction of an action of a yangian on the multiplicity space, can be obtained as  the Zhelobenko's base using a relation between the restriction problems $\mathfrak{sp}_{2n}\downarrow\mathfrak{sp}_{2n-2}$  and $\mathfrak{gl}_{n+1}\downarrow\mathfrak{gl}_{n-1}$   (see Corllary \ref{mcor}) .

Since both approaches to the construction of the base use specific
normalization  it is reasonable to cosider the vectors modulo
multiplication by a constant
\section{The Zhelobenko's construction}

\subsection{The method of $Z$-invariants of Zhelobenko}

In the present paper we use a realization of $\mathfrak{sp}_{2n}$ in which it is spanned by matrices

$$F_{i,j}=E_{i,j}-sign(i)sign(j)E_{-j,-i},\,\,\,i,j=-n,...,n,$$

The only linear relation between these matrices is the following $F_{i,j}=-sign(i)sign(j)F_{-j,-i}$, also

$$F_{i,j}\text{ is  }\begin{cases} i=j \text{ a Cartan element } ,\\ i>j \text{ a negative root element },\\
i<j  \text{ a positive root element  } \end{cases}$$

On a dense open subset of $Sp_{2n}^0$ in the group $Sp_{2n}$ we have
a Gauss decomposition

\begin{align}
\begin{split}
\label{gaude}&Sp_{2n}^0=Z^{-}DZ,\,\,\, X=\zeta \delta z,\,\,\,\\& X\in
Sp_{2n},\,\,\, \zeta\in Z^{-},\,\,\, \delta\in D,\,\,\, z\in
Z,
\end{split}
\end{align}

where  $Z^{-}$  is a subgroup of lower-triangular matrices from $Sp_{2n}$ with unit on the diagonal,
$D$ is a subgroup of diagonal matrices  from $Sp_{2n}$ ,   $Z$ is a subgroup of upper-triangular matrices
  from $Sp_{2n}$  with units on the diagonal. Let $z=(z_{i,j})$, $\delta=diag(\delta_i)$. The exists an action of  $Sp_{2n}$ on the function on  $Z$ by the following ruler.
    Let us be given a function on $Z$ of type $f(z)=f(z_{i,j})$, $i<j$.  For $X\in Sp_{2n}$ put

\begin{align*}&(Xf)(z)=\alpha(\widetilde{\delta})f(\widetilde{z}),\,\,\,
zX=\widetilde{\zeta}\widetilde{\delta}\widetilde{z}\\
&\alpha(\delta)=\delta_{-1}^{r_{-1}}...\delta_{-n}^{r_{-n}},\end{align*} and  $r_i$
 are some non-negative integers. To obtain below a realization of an
  irrep with the highest weight  $[m_{-n},...,m_{-1}]$  we must put $r_{-i}=m_{-i}-m_{-i+1}$,
  $r_{-1}=m_{-1}$. Here $m_{-i}$ is an eigenvalue on the highest vector of the operator $F_{-i,-i}$.

All function on    $Z$  form a reducible representation, to obtain
an irrep with the highest weight $[m_{-n},...,m_{-1}]$ let us do the
following. Let  $L_{i,j}$,  $i<j$  be an operator on the functions
$f(z)$ which makes a left infinitesimal shift of the function $f(z)$
to  $F_{i,j}$. Then an  irrep with the highest weight $[m_1,...,m_n]
$ is formed by functions that satisfy the {\it indicator system}
$$L_{-n,-n+1}^{r_{-n}+1}f=0,...,L_{-1,1}^{r_{-1}+1}f=0.$$

In an explicit form it it a system of PDE.

The construction of the Gelfand-Tsetlin type  base is based on the
investigation of the branching of an irrep under the restriction of
the algebra. The Zhelobenko's realization is very useful in the
investigation of this problem.

Let $V$  be an  $\mathfrak{sp}_{2n}$-irrep with the highest weight
 $[m_{-n},...,m_{-1}]$. We restrict the algebra
$\mathfrak{sp}_{2n}\downarrow \mathfrak{sp}_{2n-2}$, where
$\mathfrak{sp}_{2n-2}$ preserves coordinates $\pm n,...,\pm 2$. Then
$V$ decomposes into  the direct sum  of
$\mathfrak{sp}_{2n-2}$-irreps  $V_{\mu}$ with the highest weight
$\mu$, every such irrep occurs with the multiplicity  $m_{\mu}$:

\begin{equation}
\label{vmu} V=\bigoplus_{\mu} m_{\mu}V_{\mu}
\end{equation}

It is clear that   $v$ is a  $\mathfrak{sp}_{2n-2}$-highest vector
if and only if  $Z_{Sp_{2n-2}}v=v$, where $Z_{Sp_{2n-2}}$ is a subgroup in  $Z$ formed by elements from $Sp_{2n-2}$.  In \cite{zh1}, \cite{zh2} that
the function $f$ is a  $Z_{Sp_{2n-2}}$-invariant if and only if
$$f=f(z_{-n,-1},...,z_{-2,-1},z_{-n,1},...,z_{-1,1}).$$

\subsection{The relation between the problems of restriction $\mathfrak{gl}_{n+1}\downarrow\mathfrak{gl}_{n-1}$ and $\mathfrak{sp}_{2n}\downarrow\mathfrak{sp}_{2n-2}$}
\label{sviaz}

Consider the algebra  $\mathfrak{gl}_{n+1}$ that acts in the space with coordinates   $-n,...,-1,1$.  For this algebra also these exists a
 realization of representations on the functions on the subgroup in $GL_{n+1}$ of upper-triangular matrices  with units in the diagonal.
  An irrep is  described by an indicator system (see \cite{zh1}).

Consider the subalgebra $\mathfrak{gl}_{n-1}$, that preserve the coordinates $-n,...,-2$. Let us be given an
irrep of  $\mathfrak{gl}_{n+1}$ with the highest weight $[m_{-n},...,m_{-1},0]$ and consider the description of the branching of this representation under the restriction $\mathfrak{gl}_{n+1}\downarrow\mathfrak{gl}_{n-1}$.

It turns out that (see \cite{zh1}) $\mathfrak{gl}_{n-1}$-highest
weight vector are functions
$$f=f(z_{-n,-1},...,z_{-2,-1},z_{-n,1},...,z_{-1,1}).$$

To obtain an irrep  with the highest weight written above one must consider the function that satisfy the indicator system. It turns out that
 being restricted to the functions  $f=f(z_{-n,-1},...,z_{-2,-1},z_{-n,1},...,z_{-1,1})$ the indicator system for
   $\mathfrak{gl}_{n+1}$ with  the highest weight $[m_{-n},...,m_{-1},0]$  and for  $\mathfrak{sp}_{2n}$
   with the highest weight   $[m_{-n},...,m_{-1}]$  coincide. Thus the sets of polynomials in variables
   $z_{-n,-1},...,z_{-2,-1},z_{-n,1},...,z_{-1,1}$ that define the $\mathfrak{gl}_{n-1}$-highest
   and the $\mathfrak{sp}_{2n}$-highest vectors coincide.   Thus the restriction
    problems $\mathfrak{gl}_{n+1}\downarrow\mathfrak{gl}_{n-1}$ and $\mathfrak{sp}_{2n}\downarrow\mathfrak{sp}_{2n-2}$ are equivalent.

\subsection{The Gelfand-Tsetlin-Zhelobenko base}

There exist a base in the space of $\mathfrak{gl}_{n-1}$-highest vectors, called the Gelfand-Tsetlin base,  whose elements are index by tableaux (see \cite{zh2})

\begin{align}
\begin{split}
\label{dia}
&m_{-n,n}\,\,\,\,\,\,\,
m_{-n+1,n}\,\,\,\,\,\,\,\,\,\,\,\,\,\,...\,\,\,\,\,\,\,m_{-1,n}\,\,\,\,\,\,\, 0\\
&\,\,\,\,\,\,\,\,\,\,m'_{-n,n}\,\,\,\,\,\,\,m'_{-n+1,n}\,\,\,\,\,\,\,...\,\,\,\,\,\,\,\,\,\,\,m'_{-1,n}\\
&\,\,\,\,\,\,\,\,\,\,\,\,\,\,\,\,\,\,\,\,m_{-n,n-1}\,\,\,\,\,\,...\,\,m_{-2,n-1}\\
\end{split}
\end{align}

where $m_{-i,n}=m_{-i}$.  Denote the tableau  \eqref{dia}  as  $(m)$.

 We do not discuss a normalization since in
the present paper we consider vectors corresponding to tableaux modulo multiplication on a
constant.
% We suggest that the normalization is chosen as in \cite{1963}.

Take polynomials corresponding to these vectors. According to the
previous Section they define al possible
$\mathfrak{sp}_{2n-2}$-highest vectors. Thus we obtain a base in the
space of $\mathfrak{sp}_{2n-2}$-highest vectors indexed by
tableau \eqref{dia}.

Thus in an $\mathfrak{sp}_{2n}$-irrep we obtain a  Gelfand-Tsetlin-Zhelobenko base, whose elements are index by tableaux of type

\begin{align}
\begin{split}
\label{diap}
&m_{-n,n}\,\,\,\,\,\,\,
m_{-n+1,n}\,\,\,\,\,\,\,\,\,\,\,\,\,\,...\,\,\,\,\,\,\,m_{-1,n}\,\,\,\,\,\,\, 0\\
&\,\,\,\,\,\,\,\,\,\,m'_{-n,n}\,\,\,\,\,\,\,m'_{-n+1,n}\,\,\,\,\,\,\,...\,\,\,\,\,\,\,\,\,\,\,m'_{-1,n}\\
&\,\,\,\,\,\,\,\,\,\,\,\,\,\,\,\,\,\,\,\,m_{-n,n-1}\,\,\,\,\,\,\,...\,\,\,\,\,m_{-2,n-1}\,\,\,\,\,\,\, \,\,\,\,\,\,\, \,\,\,\,\,\,\, \,\,  0\\
&\,\,\,\,\,\,\,\,\,\,\,\,\,\,\,\,\,\,\,\,\,\,\,\,\,\,\,\,\,\,m'_{-n,n-1}\,\,\,\,\,\,...\,\,\,\,\,\,\,\,\,\,\,m'_{-2,n-1}\\
&\,\,\,\,\,\,\,\,\,\,\,\,\,\,\,\,\,\,\,\,\,\,\,\,\,\,\,\,\,\,\,\,\,\,\,\,\,\,\,\,m_{-n,n-2}\,\,\,\,\,\,\,m_{-3,n-2}\,\,\,\, 0\\
&\,\,\,\,\,\,\,\,\,\,\,\,\,\,\,\,\,\,\,\,\,\,\,\,\,\,\,\,\,\,\,\,\,\,\,\,\,\,\,\,...\\
&\,\,\,\,\,\,\,\,\,\,\,\,\,\,\,\,\,\,\,\,\,\,\,\,\,\,\,\,\,\,\,\,\,\,\,\,\,\,\,\,\,\,\,\,\,\,\,\,\,\,\,\,\,\,\,\,\,m_{-n,1}\,\,\,\,\,\,\,\,\,\,\,0\\
&\,\,\,\,\,\,\,\,\,\,\,\,\,\,\,\,\,\,\,\,\,\,\,\,\,\,\,\,\,\,\,\,\,\,\,\,\,\,\,\,\,\,\,\,\,\,\,\,\,\,\,\,\,\,\,\,\,\,\,\,\,\,\,\,\,\,\,\,\,\,\,m'_{-n,1}\\
\end{split}
\end{align}

Denote this tableau as $[m]$.

%Let  $[m]_n$ - be the first three rows of $(m)$, let  $[m]_{n-1}$  denote the rows three, four , five e.c.t..

Denote the weight component with the number $-k$ as  $weight_{-k}[m]$.

\begin{lem}
\begin{align}
\label{forw} weight_{-(n-k+1)}[m]=2\sum m'_{i,k}-\sum m_{i,k-1}-\sum
m_{i,k}
\end{align}
\end{lem}

 \subsection{A realization on the functions on the whole group}

Below we need another realization of a representation. This is a realization in the space of the function on the whole group  $Sp_{2n}$.
Onto the function  $f(g)$  the element $X\in Sp_{2n}$  acts as follows

$$(Xf)(g)=f(gX).$$

The functions that form as irrep with a given highest weight are selected by some conditions (see
\cite{zh1}). 
% For us it is important that when we restrict these selected functions  to $Z$ we obtain functions on $Z$ that form a representation selected by the indicator system  (see
%\cite{zh1}).

%as follows. They vanish under the right infinitesimal shift by an element of $Z_{-}$ and they satisfy the indicator system
% (see
%\cite{zh1}).

  Let $x_{i}^{j}$ and $\frac{\partial}{\partial x_i^j}$, $i,j=-n,...,n$  be a function of a matrix element on the space $GL_{2n}$ of all  invertible  $(2n\times 2n)$-matrices\footnote{Thus there are no relations between these functions} and a differential operator on the space of such functions.  Below we restrict the functions onto  $Sp_{2n}\subset GL_{2n}$. Thus there appear relations between $x_{i}^{j}$, but we do not   impose relations between differential operators.

Let us define the operators

$$a_{i}^j=x_i^j\cdot ,\,\,\,\, a_{i}^{+,j}=\frac{\partial}{\partial x_i^j},\,\,\,i,j=-n,...,n$$

Then $$F_{i,j}=\sum_{k} \Big(
a_{i}^{k}a^{+,k}_{j}-sing(i)sign(j)a_{-j}^{k}a^{+,k}_{-i}\Big).$$
Also define
$$a_{i_1,...,i_k}:=det(a_i^j)_{i=i_1,...,i_k}^{j=-n,...,-n+k-1}$$
One can easily check that the vector
\begin{equation}
\label{stv} \prod_{k=-n}^{-1} (a_{-n,...,-k})^{m_{-k}-m_{-k+1}}1,
\end{equation}

where  $1$ is a constant function,
is the highest vector.

 One can realize  representation $\mathfrak{gl}_{n+1} $ in the space of
the functions on   $GL_{n+1}$ and also define  $a_{i}^j$ and
$a_{i}^{+,j}$, but for $i,j=-n,...,-1,1$. The generators  $E_{i,j}$
act by formulas
$$E_{i,j}=\sum_k a_{i}^{k}a^{+,k}_{j}.$$
The formula \eqref{stv} defines also  a $\mathfrak{gl}_{n+1}$-highest vector of weight  $[m_{-n},...,m_{-1},0]$.

\section{The Mikelsson-Zhelobenko algebra}

\subsection{The definition}
\label{defmz}
The definition of the Mikelsson-Zhelobenko algebra is given in \cite{zh3} (see also \cite{M},\cite{zh4}). %  В настоящей работе оно в полной общности не приводится.

   Let us be given a reductive Lie algebra  $\mathfrak{g}$ and it's  reductive subalgebra  $\mathfrak{k}$.  Consider the triangular decomposition
$$
\mathfrak{k}=\mathfrak{k}_{-}\oplus\mathfrak{k}_0\oplus\mathfrak{k}_+,
$$

Denote the set of positive roots for the algebra  $\mathfrak{k}$ as  $\Delta_+$, the element of  $\mathfrak{k}_+$
 corresponding to the root $\alpha\in \Delta_+$ denote as  $e_{\alpha}$. There exist also elements  $e_{-\alpha}\in\mathfrak{k}_{-}$.

Let   $R$ -be a field of fractions of   $U(\mathfrak{k}_0)$, put
$
U'(\mathfrak{g})=U(\mathfrak{g})\otimes_{U(\mathfrak{k}_0)}R.
$
Let  $J'=U'(\mathfrak{g})\mathfrak{k}_+$.

Define $$M(\mathfrak{g},\mathfrak{k}):=U'(\mathfrak{g})/J'.$$ {\it The the Mikelsson-Zhelobenko algebra} is defined as follows

$$
Z(\mathfrak{g},\mathfrak{k}):=\{v\in M(\mathfrak{g},\mathfrak{k}):\,\,\,\,\mathfrak{k}^+v=0\}.
$$

One can describe the generator of this algebra. For $\alpha\in\Delta_+$  put

$$
p_{\alpha}=1+\sum_{k=1}^{\infty}e_{-\alpha}^ke_{\alpha}\frac{(-1)^k}{k!(h_{\alpha}+\rho(h_{\alpha})+1)...(h_{\alpha}+\rho(h_{\alpha})+k)},
$$

where $h_{\alpha}=[e_{\alpha},e_{-\alpha}]$, $\rho$  is a half-sum of positive roots. Choose a linear order  $\alpha_1<...<\alpha_m$ of roots from  $\Delta_+$ such that if one root is a sum of two others than it stands between them.

Put

$$
p_{\mathfrak{g},\mathfrak{k}}:=p_{\alpha_1}...p_{\alpha_m}.
$$
\begin{prop}
The operation of multiplication onto  $p_{\mathfrak{g},\mathfrak{k}}$ on the left side  projects  $M(\mathfrak{g},\mathfrak{k})$  to  $Z(\mathfrak{g},\mathfrak{k})$.
\end{prop}

\begin{prop}
Let as a linear space $\mathfrak{g}=\mathfrak{k}\oplus<E_1,...,E_k>$, then  $Z(\mathfrak{g},\mathfrak{k})$ is generated by elements $e_i:=p_{\mathfrak{g},\mathfrak{k}}E_i$
\end{prop}

Let us be given a representation  $V$ of the Lie algebra $\mathfrak{g}$.  Consider the subspace $V^+$ that consists of  $\mathfrak{k}$-highest vectors.
 \begin{prop}
 For $v\in V$ one has
\begin{align*}p_{\mathfrak{g},\mathfrak{k}}v=
\begin{cases}
v,\,\,\, \text{  if }v \text{  is a $\mathfrak{k}$-highest vector },\\
0,\,\,\, \text{ if }v=e_{-}w,\text{ where  }e_{-}\in \mathfrak{k}
\text{  negative root element }
\end{cases}
\end{align*}
%  $pv=0$, if  $v$ is not a  $\mathfrak{k}$-highest vector and   $pv=v$, if  $v$ is a   $\mathfrak{k}$-highest vector.
 \end{prop}

  \begin{prop}
The action  of   $\mathfrak{g}$ on $V$ induces an action $Z(\mathfrak{g},\mathfrak{k})$ on $V^+$.
 \end{prop}

%Note the following fact (see the Section  \ref{yafo}).

% В разделе  \ref{yafo} был определён экстремальный проектор
% проектор $p_{\mathfrak{g},\mathfrak{k}}$, действующий в представлениях  $\mathfrak{g}$. Определённый в настоящем разделе оператор  $p$ действует в представлениях именно так, как описано в разделе \ref{yafo}.

\subsection{ The algebras used in the paper}

\subsubsection{$Z^{Zh}(\mathfrak{gl}_{n+1},\mathfrak{gl}_{n-1})$}

Let the algebra $\mathfrak{gl}_{n+1}$act in the space with
coordinates $-n,...-2,-1,1$, and let $\mathfrak{gl}_{n-1}$ be a
subalgebra that preserves the coordinates $-n,...-2$.

Thus  $Z^{Zh}(\mathfrak{gl}_{n+1},\mathfrak{gl}_{n-1})$ is an
associative algebra with generators $e_{\pm 1,i}$, $e_{i,\pm 1}$,
$i=-n,...,-2$, and also  $e^{\pm 1}_{j,j}$,
$j=-n,...,-2$.  Here $ e_{i,j}:=p_{\mathfrak{gl}_{n+1},\mathfrak{gl}_{n-1}}E_{i,j}$.

This algebra has the following representations  For an
arbitrary representation $\mathfrak{gl}_{n+1}$ the algebra
$Z^{Zh}(\mathfrak{gl}_{n+1},\mathfrak{gl}_{n-1})$ acts on the space
of $\mathfrak{gl}_{n-1}$-highest vectors by the ruler

\begin{align*}
& e_{i,j}\mapsto p_{\mathfrak{gl}_{n+1},\mathfrak{gl}_{n-1}}E_{i,j},  \,\,\,  e_{i,i}^{\pm 1}\mapsto E^{\pm 1}_{i,i},
\end{align*}

Note that the action of $E_{i,i}^{-1}$ is well defined.

Define also elements $e_{-1,1}=p_{\mathfrak{gl}_{n+1},\mathfrak{gl}_{n-1}}E_{-1,1}$, $e_{1,-1}p_{\mathfrak{gl}_{n+1},\mathfrak{gl}_{n-1}}E_{1,-1}$.
\subsubsection{$Z^{Zh}(\mathfrak{sp}_{2n},\mathfrak{sp}_{2n-2})$ }

Chose a subalgebra $\mathfrak{sp}_{2n-2}$  that preserves the
coordinates  $\{-n,...,-2,2,...,n\}$.  The algebra
$Z^{Zh}(\mathfrak{sp}_{2n},\mathfrak{sp}_{2n-2})$ is an associative
algebra generated by  $\zeta_{\pm 1,i}$, $\zeta_{i,\pm 1}$,
$i=-n,...,-2$ and also  , $\zeta^{\pm
1}_{j,j}$,  $j=-n,...,-2,2,...,n$. Here  $\zeta_{i,j}= p_{\mathfrak{sp}_{2n},\mathfrak{sp}_{2n-2}}F_{i,j}$. This algebra has the following
representations. For an arbitrary representation  of
$\mathfrak{sp}_{2n}$   the algebra
$Z^{Zh}(\mathfrak{sp}_{2n},\mathfrak{sp}_{2n-2})$ acts on the space
of $\mathfrak{sp}_{2n-2}$-highest vectors by the ruler

\begin{align*}
& \zeta_{i,j}\mapsto p_{\mathfrak{sp}_{2n},\mathfrak{sp}_{2n-2}}F_{i,j},  \,\,\,  F_{i,i}^{\pm 1}\mapsto F^{\pm 1}_{i,i},
\end{align*}

One can also add elements
$\zeta_{i,j}=sign(i)sign(j)\zeta_{-j,-i}$.

Also put $\zeta_{-1,1}= p_{\mathfrak{sp}_{2n},\mathfrak{sp}_{2n-2}}F_{-1,1}$, $\zeta_{-1,1}= p_{\mathfrak{sp}_{2n},\mathfrak{sp}_{2n-2}}F_{1,-1}$.

\subsubsection{$Z^{M}(\mathfrak{gl}_{n+1},\mathfrak{gl}_{n-1})$}

Let the algebra $\mathfrak{gl}_{n+1}$act in the space with
coordinates $-n,-n+1...,-1,n$, and let $\mathfrak{gl}_{n-1}$ be a
subalgebra that preserves the coordinates $-n+1,...-1$.

Thus  $Z^{M}(\mathfrak{gl}_{n+1},\mathfrak{gl}_{n-1})$ is an
associative algebra with generators $\varepsilon_{\pm n,i}$,
$\varepsilon_{i,\pm n}$, $i=-n+1,...,-1$, and also $\varepsilon_{-n,n}$,
$\varepsilon_{n,-n}$, $\varepsilon^{\pm 1}_{j,j}$, $j=-n+1,...,-1$. Here  $\varepsilon_{i,j}= p_{\mathfrak{gl}_{n+1},\mathfrak{gl}_{n-1}}E_{i,j}$.

Also put $\varepsilon_{-n,n}= p_{\mathfrak{gl}_{n+1},\mathfrak{gl}_{n-1}}E_{-n,n}$, $\varepsilon_{n,-n}= p_{\mathfrak{gl}_{n+1},\mathfrak{gl}_{n-1}}E_{n,-n}$

\subsubsection{$Z^{M}(\mathfrak{sp}_{2n},\mathfrak{sp}_{2n-2})$ }

Chose a subalgebra  $\mathfrak{sp}_{2n-2}$ that preserves the coordinates $\{-n+1,...,n-1\}$.
The algebra $Z^{М}(\mathfrak{sp}_{2n},\mathfrak{sp}_{2n-2})$ is an associative subalgebra generated by elements
 $\eta_{\pm n,i}$, $\eta_{i,\pm n}$, $i=-n+1,...,-1$, and also  $\eta^{\pm 1}_{j,j}$,  $j=-n+1,...,n-1$, where  $\eta_{i,j}= p_{\mathfrak{sp}_{2n},\mathfrak{sp}_{2n-2}}F_{i,j}$.

 This algebra has the following representations. For an arbitrary representation  of $\mathfrak{sp}_{2n}$
   the algebra $Z^{Zh}(\mathfrak{sp}_{2n},\mathfrak{sp}_{2n-2})$ acts on the space of $\mathfrak{sp}_{2n-2}$-highest vectors by the ruler

\begin{align*}
& \eta_{i,j}\mapsto p_{\mathfrak{sp}_{2n},\mathfrak{sp}_{2n-2}}F_{i,j},  \,\,\,  F_{i,i}^{\pm 1}\mapsto F^{\pm 1}_{i,i},
\end{align*}
One can also add elements  $\eta_{i,j}=sign(i)sign(j)\eta_{-j,-i}$.

Define also elements  $\eta_{-n,n}=p_{\mathfrak{sp}_{2n},\mathfrak{sp}_{2n-2}}F_{-n,n}$, $\eta_{n,-n} =p_{\mathfrak{sp}_{2n},\mathfrak{sp}_{2n-2}}F_{n,-n}$.

For $i=1,...,n$ define $\rho_{-i}=-\rho_i=i$, $f_{i}=F_{i,i}+\rho_i$ and  $f_{-i}=-\rho_i$. Let $g_i=f_{i}+\frac{1}{2}$ for all  $i$.

Define the operators \begin{equation}\label{ck}\check{\eta}_{i,\pm n}=\eta_{i,\pm n}(f_i-f_{i-1})...(f_i-{-n+1}).\end{equation}

Define elements $Z_{n,-n}(u)$ by formula

\begin{align}
\begin{split}
\label{znn} &Z_{n,-n}(u)=-F_{n,-n}\prod_{i=-n+1}^{n-1}(u+g_i)+
\sum_{i=-n+1}^{n-1}\check{\eta}_{n,i}\check{\eta}_{i,n}\prod_{j=-n+1,j\neq
}^{n-1}\frac{u+g_j}{g_i-g_j}
\end{split}
\end{align}

\subsection{The action of the Mikelsson-Zhelobenko algebra onto the determinants}

In the Corollary \ref{slea} below we show that in the restriction problems $\mathfrak{gl}_{n+1}\downarrow\mathfrak{gl}_{n-1}$ and $\mathfrak{sp}_{2n}\downarrow\mathfrak{sp}_{2n-2}$ in the Zhelobenko's approach
the corresponding highest vector is a polynomial in  $a_{-n,...,-k}$, $a_{-n,...-k-1,-1}$, $a_{-n,...,-k-1,1}$, $a_{-n,...-k-2,-1,1}$.

And in the Molev's approach - a polynomial in $a_{-n,...,-k}$, $a_{-n,...-k-1,n}$, $a_{-n+1,...,-k+1}$, $a_{-n+1,...-k,n}$.

Let us find some relations between operators, by means of which the generators of the Mikelsson-Zhelobenko algebra act on the determinants.

\subsubsection{$Z^{Zh}(\mathfrak{gl}_{n+1},\mathfrak{gl}_{n-1})$}

The action of the generators onto determinants is described by the following diagram

\begin{displaymath}
 \xymatrix{a_{-n} \ar[r]^{e_{-1,-n}}
 \ar[d]^{e_{1,-n}}
   & a_{-1} \ar@<1ex>[l]^{e_{-n,-1}}   \ar@<1ex>[ld] \\
a_{1}  \ar@<1ex>[u]^{e_{-n,1}} \ar[ur]  }
\end{displaymath}

An arrow that goes by diagonal up and right  is $e_{1,-1}$,  and an arrow that goes left and down is  $e_{-1,1}$.

For $1<k<n$ one has

\begin{displaymath}
 \xymatrix{a_{-n,\,\,\,\,.\,.\,.\,\,\,\,,-k} \ar[r]^{e_{-1,-k}}
 \ar[d]^{e_{1,-k}}
   & a_{-n,\,\,\,\,.\,.\,.\,\,\,\,,-k-1,-1} \ar@<1ex>[l]^{e_{-k,-1}}  \ar [d]^{-e_{1,-k-1}}  \ar@<1ex>[ld] \\
a_{-n,\,\,\,\,.\,.\,.\,\,\,\,,-k-1,1} \ar[r]^{e_{1,-k-1}}
\ar@<1ex>[u]^{e_{-k,1}} \ar[ur] &
a_{-n,\,\,\,\,.\,.\,.\,\,\,\,-k-2,-1,1} \ar@<1ex>[l]^{e_{-k-1,-1}}
\ar@<1ex>[u]^{-e_{-k-1,1}} }
\end{displaymath}

An arrow that goes by diagonal up and right  is $e_{1,-1}$,  and an arrow that goes left and down is  $e_{-1,1}$.

For $k=1$

\begin{displaymath}
 \xymatrix{
   & a_{-n,\,\,\,\,.\,.\,.\,\,\,\,,-2,-1}   \ar [d]^{-e_{1,-2}}  \ar@<1ex>[ld]^{e_{1,-1}} \\
a_{-n,\,\,\,\,.\,.\,.\,\,\,\,,-2,1} \ar[r]^{e_{1,-2}}
 \ar[ur]^{e_{-1,1}} &
a_{-n,\,\,\,\,.\,.\,.\,\,\,\,-3,-1,1} \ar@<1ex>[l]^{e_{-2,-1}}
\ar@<1ex>[u]^{-e_{-2,1}} }
\end{displaymath}

\subsubsection{ $Z^{Zh}(\mathfrak{sp}_{2n},\mathfrak{sp}_{2n-2})$}

The action of the generators to the determinant is described by the following diagrams

\begin{displaymath}
 \xymatrix{a_{-n} \ar[r]^{\zeta_{-1,-n}}
 \ar[d]^{\zeta_{1,-n}}
   & a_{-1} \ar@<1ex>[l]^{\zeta_{-n,-1}}   \ar@<1ex>[ld] \\
a_{1}  \ar@<1ex>[u]^{\zeta_{-n,1}} \ar[ur]  }
\end{displaymath}

An arrow that goes by diagonal up and right  is
$\zeta_{1,-1}$,  and an arrow that goes left and down is
$\zeta_{-1,1}$.

For  $1<k<n$

\begin{displaymath}
 \xymatrix{a_{-n,\,\,\,\,.\,.\,.\,\,\,\,,-k} \ar[r]^{\zeta_{-1,-k}}
 \ar[d]^{\zeta_{1,-k}}
   & a_{-n,\,\,\,\,.\,.\,.\,\,\,\,,-k-1,-1} \ar@<1ex>[l]^{\zeta_{-k,-1}}  \ar [d]^{-\zeta_{1,-k-1}}  \ar@<1ex>[ld] \\
a_{-n,\,\,\,\,.\,.\,.\,\,\,\,,-k-1,1} \ar[r]^{\zeta_{-1,-k-1}}
\ar@<1ex>[u]^{\zeta_{-k,1}} \ar[ur] &
a_{-n,\,\,\,\,.\,.\,.\,\,\,\,-k-2,-1,1}
\ar@<1ex>[l]^{\zeta_{-k-1,-1}} \ar@<1ex>[u]^{-\zeta_{-k-1,1}} }
\end{displaymath}

An arrow that goes by diagonal up and right  is
$\zeta_{1,-1}$,  and an arrow that goes left and down is
$\zeta_{-1,1}$.

For $k=1$ we have  the diagram

\begin{displaymath}
 \xymatrix{
   & a_{-n,\,\,\,\,.\,.\,.\,\,\,\,,-2,-1}   \ar@{=>}[d]^{-\zeta_{1,-2}}  \ar@<1ex>[ld]^{} \\
a_{-n,\,\,\,\,.\,.\,.\,\,\,\,,-2,1}
\ar@{=>}[r]^{\zeta_{1,-2}}
 \ar[ur]^{} &
a_{-n,\,\,\,\,.\,.\,.\,\,\,\,-3,-1,1}
\ar@<1ex>[l]^{\zeta_{-2,-1}}
\ar@<1ex>[u]^{-\zeta_{-2,1}} }
\end{displaymath}

An arrow that goes by diagonal up and right  is
$\zeta_{1,-1}$,  and an arrow that goes left and down is
$\zeta_{-1,1}$.
A double arrow means that one determinant is mapped to another multiplied by
 $2$.

To obtain this diagram we use the following fact. We have

\begin{align*}
&\zeta_{1,-2}a_{-n,...,-2,-1}=a_{-n,...,-3,1,-1}+a_{-n,....,-3,-2,2},\\
&\zeta_{-1,-2}a_{-n,...,-2,1}=a_{-n,...,-3,-1,1}-a_{-n,....,-3,-2,2},\\
\end{align*}

But as it is shown below due to the relations between matrix elements in the group
$Sp_{2n}$ we have $a_{-n,....,-3,-2,2}=-a_{-n,...,-3,-1,1}$.  This is the reason of the appearance of twin arrows.

\begin{lem}
 For functions on $Sp_{2n}$ one has a relation $a_{-n,....,-3,-2,2}=-a_{-n,...,-3,-1,1}$
\end{lem}
\proof

First of all we prove the equality for the group   $Sp_{4}$. It can be done as follows. Take a matrix  $X\in Sp_{4}^0$ and consider it's decomposition  \eqref{gaude}.
The matrices  $\zeta$, $\delta$ and  $z$ can be represented as exponents of Lie algebra elements $\zeta=e^{A}$, $\delta=e^B$, $z=e^C$. The matrix  $A$is a linear combination of  $F_{i,j}$, $i>j$,  $B$ is a linear combination of   $F_{i,i}$, $C$  is a linear combination of    $F_{i,j}$, $i<j$. Using this fact we obtain  a parametrization of matrices $A$, $B$, $C$, then after taking an exponent a parametrization of matrices $\zeta$, $\delta$ и  $z$, and finally a presentation of an arbitrary matrix  $X\in Sp_{4}^0$. Then we can check the equality $a_{-n,....,-3,-2,2}=-a_{-n,...,-3,-1,1}$ by direct computations \footnote{Of course it is better to do it using a computer}. Since   $Sp_{4}^0$ is dense in    $Sp_{4}$ the equality holds everywhere on  $Sp_{4}$.

Consider the case of an arbitrary group  $Sp_{2n}$.   Let  $\mathcal{F}_{i,j}(\alpha)$  be a matrix with units on the diagonal, with   $\alpha$ on the place  $(i,j)$  and with  $-alpha$ on the place $(-j,-i)$.  The matrices  $\mathcal{F}_{i,j}(\alpha)$ belong to $Sp_{2n}$. A multiplication by this matrix on the right is equivalent to  doing an elementary transformation of the matrix $X$: we add to the  $j$-th row the $i$-th row with a coefficient  $\alpha$, and we add simultaneously to the $-i$-th row the $-j$-the row with a coefficient $-\alpha$.  A multiplication by $\mathcal{F}_{i,j}(\alpha)$  on the left is equivalent to an analogous transformation of rows.
% Отметим, что выполнение такого преобразования со строками  $-т,-4,...,-n$  не меняет обоих определителей  $a_{-n,....,-3,-2,2}$ и $a_{-n,...,-3,-1,1}$.

Multiplying by  $\mathcal{F}_{i,j}(\alpha)$ on the left and on the right we can transform   $X$ preserving  $a_{-n,....,-3,-2,2}$ and $a_{-n,...,-3,-1,1}$ to the following form

\begin{equation}
\begin{pmatrix}
x_{-n,-n}&...&0& 0 & 0& 0 &0& 0 &...& x_{-n,n}\\
...\\
0&...&x_{-3,-3}& 0 & 0 &0 &0& 0& ...&0\\
0&...&0& x_{-2,-2} & x_{-2,-1} & x_{-2,1} & x_{-2,2}&0&...&0\\
0&...&0& x_{-1,-2} & x_{-1,-1} & x_{-1,1} & x_{-1,2}&0&...&0\\
0&...&0& x_{1,-2} & x_{1,-1} & x_{1,1} & x_{1,2}&0&...&0\\
0&...&0& x_{2,-2} & x_{2,-1} & x_{2,1} & x_{2,2}&0&...&0\\
0&...&x_{3,-3}& 0 & 0 & 0 &0&x_{3,3}&...&0\\
...\\
x_{n,-n}&...&0& 0 & 0 & 0 &0&0&...&x_{n,n}\\
\end{pmatrix}
\end{equation}

Since matrices  $\mathcal{F}_{i,j}(\alpha)$  belong to $Sp_{2n}$, this matrix belongs to  $Sp_{2n}$.  We have

\begin{align*}
&a_{-n,....,-3,-2,2}=x_{-n,-n}...x_{-3,-3}\cdot det\begin{pmatrix}  x_{-2,-2} & x_{-2,2}\\ x_{-1,-2} &x_{-1,2}\end{pmatrix},\\
&a_{-n,....,-3,-1,1}=x_{-n,-n}...x_{-3,-3}\cdot det\begin{pmatrix} x_{-2,-1} & x_{-2,1}\\ x_{-1,-1} &x_{-1,1}\end{pmatrix}.
\end{align*}

The submatrix

\begin{equation}
\begin{pmatrix}

x_{-2,-2} & x_{-2,-1} & x_{-2,1} & x_{-2,2}\\
 x_{-1,-2} & x_{-1,-1} & x_{-1,1} & x_{-1,2}\\
 x_{1,-2} & x_{1,-1} & x_{1,1} & x_{1,2}\\
 x_{2,-2} & x_{2,-1} & x_{2,1} & x_{2,2}\\
\end{pmatrix}
\end{equation}

 belongs to $Sp_{4}$. Using the equality  $a_{-2,2}=-a_{-1,1}$ for  $Sp_{4}$ we obtain that   $a_{-n,....,-3,-2,2}$ и $a_{-n,...,-3,-1,1}$ в $Sp_{2n}$.

\endproof

\subsubsection{$Z^{M}(\mathfrak{gl}_{n+1},\mathfrak{gl}_{n-1})$}

The action of the generators to the determinant is described by the
following diagrams

\begin{displaymath}
 \xymatrix{a_{-n} \ar[r]^{\varepsilon_{n,-n}}
 \ar[d]^{\varepsilon_{-n+1,-n}}
   & a_{n} \ar@<1ex>[l]^{\varepsilon_{-n,n}}  \ar [ld]   \\
a_{-n+1} \ar@<1ex>[u]^{\varepsilon_{-n,-n+1}} \ar@<1ex>[ur] &
    }
\end{displaymath}

An arrow that goes by diagonal up and right  is $\eta_{n,-n+1}$,
and an arrow that goes left and down is $\eta_{-n+1,n}$.

For $1<k<n$ we have

\begin{displaymath}
 \xymatrix{a_{-n,\,\,\,\,.\,.\,.\,\,\,\,,-k} \ar[r]^{\varepsilon_{n,-k}}
 \ar[d]^{\varepsilon_{-k+1,-n}} \ar@<1ex>[rd]
   & a_{-n,\,\,\,\,.\,.\,.\,\,\,\,,-k-1,n} \ar@<1ex>[l]^{\varepsilon_{-k,n}}  \ar [d]^{\varepsilon_{-k,-n}}   \\
(-1)^{n-k-1}a_{-n+1,\,\,\,\,.\,.\,.\,\,\,\,,-k,-k+1}
\ar[r]^{\varepsilon_{n,-k+1}} \ar@<1ex>[u]^{\varepsilon_{-n,-k+1}} &
(-1)^{n-k-1}a_{-n+1,\,\,\,\,.\,.\,.\,\,\,\,-k,n}
\ar@<1ex>[l]^{\varepsilon_{-k+1,n}} \ar@<1ex>[u]^{\varepsilon_{-n,-k}}
\ar[ul] }
\end{displaymath}

An arrow that goes by diagonal up and right  is
$(-1)^{n-k-1}\varepsilon_{n,-n}$,   and an arrow that goes left and down is
$(-1)^{n-k-1}\varepsilon_{-n,n}$.

Now let $k=1$.

\begin{displaymath}
 \xymatrix{a_{-n,\,\,\,\,.\,.\,.\,\,\,\,,-1}
 \ar@{->}[d]^{\varepsilon_{1,-n}} \ar@<1ex>[rd]_{}
   &     \\
(-1)^{n-1}a_{-n+1,\,\,\,\,.\,.\,.\,\,\,\,,-1,1}
\ar@<1ex>[r]^{\varepsilon_{n,1}} \ar@<1ex>[u]^{\varepsilon_{-n,1}} &
(-1)^{n-1}a_{-n+1,\,\,\,\,.\,.\,.\,\,\,\,-1,n} \ar@{->}[l]^{\varepsilon_{1,n}}
  \ar[ul]_{}   }
\end{displaymath}

An arrow that goes by diagonal up and right  is
$(-1)^{n-1}\varepsilon_{n,-n}$,   and an arrow that goes left and down is
$(-1)^{n-1}\varepsilon_{-n,n}$.

\subsubsection{ $Z^{M}(\mathfrak{sp}_{2n},\mathfrak{sp}_{2n-2})$}

The action of the generators to the determinant is described by the following diagrams

\begin{displaymath}
 \xymatrix{a_{-n} \ar[r]^{\eta_{n,-n}}
 \ar[d]^{\eta_{-n+1,-n}}
   & a_{n} \ar@<1ex>[l]^{\eta_{-n,n}}  \ar [ld]   \\
a_{-n+1} \ar@<1ex>[u]^{\eta_{-n,-n+1}} \ar@<1ex>[ur] &
    }
\end{displaymath}

An arrow that goes by diagonal up and right  is
$\eta_{n,-n+1}$,    and an arrow that goes left and down is
$\eta_{-n+1,n}$.

For $1<k<n$ we have

\begin{displaymath}
 \xymatrix{a_{-n,\,\,\,\,.\,.\,.\,\,\,\,,-k} \ar[r]^{\eta_{n,-k}}
 \ar[d]^{\eta_{-k+1,-n}} \ar@<1ex>[rd]
   & a_{-n,\,\,\,\,.\,.\,.\,\,\,\,,-k-1,n} \ar@<1ex>[l]^{\eta_{-k,n}}  \ar [d]^{\eta_{-k,-n}}   \\
(-1)^{n-k-1}a_{-n+1,\,\,\,\,.\,.\,.\,\,\,\,,-k,-k+1} \ar[r]^{\eta_{n,-k+1}}
\ar@<1ex>[u]^{\eta_{-n,-k+1}} & (-1)^{n-k-1}a_{-n+1,\,\,\,\,.\,.\,.\,\,\,\,-k,n}
\ar@<1ex>[l]^{\eta_{-k+1,n}} \ar@<1ex>[u]^{\eta_{-n,-k}}  \ar[ul] }
\end{displaymath}

An arrow that goes by diagonal up and right  is $(-1)^{n-k-1}\eta_{n,-n}$,   and an arrow that goes left and down is $(-1)^{n-k-1}\eta_{-n,n}$.

Now let $k=1$.

\begin{displaymath}
 \xymatrix{a_{-n,\,\,\,\,.\,.\,.\,\,\,\,,-1}
 \ar@{=>}[d]^{\eta_{1,-n}} \ar@<1ex>[rd]_{}
   &     \\
(-1)^{n-1}a_{-n+1,\,\,\,\,.\,.\,.\,\,\,\,,-1,1}
\ar@<1ex>[r]^{\eta_{n,1}} \ar@<1ex>[u]^{\eta_{-n,1}}
& (-1)^{n-1}a_{-n+1,\,\,\,\,.\,.\,.\,\,\,\,-1,n}
\ar@{=>}[l]^{\eta_{1,n}}
  \ar[ul]_{}   }
\end{displaymath}

An arrow that goes by diagonal up and right  is $ (-1)^{n-1}\eta_{n,-n}$,   and an arrow that goes left and down is $ (-1)^{n-1}\eta_{-n,n}$.

\subsubsection{The action of  $Z_{n,-n}(u)$}

Let us discuss the action of $Z_{n,-n}(u)$ onto determinants.

\begin{lem}\label{znnac}
The operator $Z_{n,-n}(u)$ acts onto determinants by ruler

\begin{align*}
&a_{-n,...,-i}\mapsto  const\cdot a_{-n+1,...,-i,n},\\
&\text{other determinants }\mapsto 0,
\end{align*}
and onto a product of determinants it acts by Leibnitz ruler. The
constant depends on $i$ and a formal variable $u$.
\end{lem}
\proof The first summand in \eqref{znn} act as it is described in
the statement of Lemma. Now consider the action of summands of type
$\check{\eta}_{n,i}\check{\eta}_{i,-n}$ for $i=1,..,n-1$.  From the
diagrams above one can make a conclusion that onto a determinant
this composition acts by ruler
\begin{align*}
&a_{-n,...,-i}\mapsto  (-1)^{n-i}a_{-n+1,...,-i,n},\\
&\text{other determinants }\mapsto 0,
\end{align*}

and onto a product of determinants it acts by Leibnitz ruler.

\endproof

\subsection{A simplified Mikelsson-Zhelobenko algebra}

Introduce three simplified Mikelsson-Zhelobenko algebras
$Z^{Zh}_s(\mathfrak{gl}_{n+1},\mathfrak{gl}_{n-1})$,
$Z^{Zh}_s(\mathfrak{sp}_{2n},\mathfrak{sp}_{2n-2})$,
$Z^{M}_s(\mathfrak{gl}_{n+1},\mathfrak{gl}_{n-1})$,
$Z^{M}_s(\mathfrak{sp}_{2n},\mathfrak{sp}_{2n-2})$.

The first two are algebras of operators that act on the space of
polynomials in  determinants $a_{-n,...,-k}$, $a_{-n,...-k-1,-1}$,
$a_{-n,...,-k-1,1}$, $a_{-n,...-k-2,-1,1}$ and satisfy the Leibnitz
ruler and commutation relations presented in the diagrams.

The third  and the fourth algebras are constructed analogously but
they are the algebras of operators that act on polynomials in
variables $a_{-n,...,-k}$, $a_{-n,...-k-1,n}$, $a_{-n+1,...,-k+1}$,
$a_{-n+1,...-k,n}$.

The  simplified Mikelsson-Zhelobenko algebra is a factor if the
ordinary  Mikelsson-Zhelobenko algebra. The action of the ordinary
Mikelsson-Zhelobenko algebra on the space of higest vectors factors
through the homomorphism into the simplified Mikelsson-Zhelobenko
algebra.

%\subsection{The comparison of actions}

Let us compare the actions of
$Z^{Zh}(\mathfrak{gl}_{n+1},\mathfrak{gl}_{n-1})$ and
$Z^{Zh}(\mathfrak{sp}_{2n},\mathfrak{sp}_{2n-2})$,
$Z^{M}(\mathfrak{gl}_{n+1},\mathfrak{gl}_{n-1})$ and
$Z^{M}(\mathfrak{sp}_{2n},\mathfrak{sp}_{2n-2})$. We come to the
following conclusion. Under the correspondence

\begin{align}
\begin{split}
\label{zhz} & e_{\pm 1,-2}\leftrightarrow \frac{1}{2}\zeta_{\pm
1,-2}\\
 & e_{i,j}\leftrightarrow \zeta_{i,j}\text{ in other cases }
\end{split}
\end{align}

\begin{align}
\begin{split}
\label{zm} & \varepsilon_{1,\pm n}\leftrightarrow
\frac{1}{2}\eta_{1,\pm n }\\
 & \varepsilon_{i,j}\leftrightarrow \eta_{i,j}\text{ in other cases }
\end{split}
\end{align}

the results of the actions on the determinants coincide.

\begin{lem}
\label{soo} The correspondences  \eqref{zhz}, \eqref{zm} give
isomorphisms of $Z^{Zh}_s(\mathfrak{gl}_{n+1},\mathfrak{gl}_{n-1})$
and $Z^{Zh}_s(\mathfrak{sp}_{2n},\mathfrak{sp}_{2n-2})$,
$Z^{M}_s(\mathfrak{gl}_{n+1},\mathfrak{gl}_{n-1})$ and
$Z^{M}_s(\mathfrak{sp}_{2n},\mathfrak{sp}_{2n-2})$. These isomorphism commute modulo multiplication onto a constant with the action on the space pf  $\mathfrak{sp}_{2n-2}$-highest vectors.

\end{lem}

\begin{lem} The
$\mathfrak{sp}_{2n-2}$-highest vector in the Zhelobenko's  approach
can be determined by formula \eqref{vzh1},
$\mathfrak{sp}_{2n-2}$--highest vector in the Molev's approach can
be determined by formula \eqref{vm1}
\end{lem}
\section{Explicit formulas for $\mathfrak{sp}_{2n-2}$-highest vectors}
\subsection{Zhelobenko's approach}

Let us give explicit formulas for $\mathfrak{sp}_{2n-2}$-highest vectors
in the Zhelobenko's approach.

\begin{align}
\begin{split}
\label{stvzhe} v=&\prod_{i=2}^{n}
e_{-1,-i}^{m'_{-i,n}-m_{-i,n-1}}\prod_{i=1}^{n}
e_{1,-i}^{m_{-i,n}-m'_{-i,n}}v_0,
\end{split}
\end{align}
where  $v_0$ is the highest weight vector.

Using lemma \ref{soo} we come to the following conclusion
\begin{lem} The
$\mathfrak{sp}_{2n-2}$-с-highest vectors
in the Zhelobenko's approach can be given by formula

\begin{align}
\begin{split}
\label{stvzhze} v=&\prod_{i=2}^{n}
\zeta_{-1,-i}^{m'_{-i,n}-m_{-i,n-1}}\prod_{i=1}^{n}
\zeta_{1,-i}^{m_{-i,n}-m'_{-i,n}}v_0,
\end{split}
\end{align}

\end{lem}

\subsection{Molev's approach}

Let us give formulas for
$\mathfrak{sp}_{2n-2}$-highest vectors in the Molev's  approach (see
\cite{M}, \cite{M1}),  corresponding to the tableau \eqref{dia}.

\begin{align}
\begin{split}
\label{stvm}
v=&\prod_{i=1}^{n-1} \check{\eta}_{-n,-i}^{m'_{-i,n}-m_{-i,n-1}}
 \check{\eta}_{n,-i}^{m_{-i,n}-m'_{-i,n}}\prod_{k=m'_{-n,n}+\rho_{n}+\frac{1}{2}}^{m_{-n,n}+\rho_{n}+\frac{1}{2}}Z_{n,-n}(k)v_0
\end{split}
\end{align}

Below we give a simplification of Molev's
formulas.

In \cite{M} there were obtained formulas for the action of
$F_{n,-n}$ and $F_{-n,n}$.
According to them after the action of the first operator onto a vector corresponding to a tableau we obtain a linear combination of vectors that are obtained by the change $m'_{-i,n}\mapsto m'_{-i,n}-1$. And after the second -   a linear combination of vectors that are obtained by the change  $m'_{-i,n}\mapsto m'_{-i,n}+1$.

\begin{lem}
The formula \eqref{stvm}   can be simplified as follows.

\begin{equation}
\label{stvm1}
v=const\prod_{i=1}^{n-1} \eta_{-n,-i}^{m'_{i,n}-m_{i,n-1}}
\prod_{i=1}^{n} \eta_{n,-i}^{m_{-i,n}-m'_{-i,n}}
\end{equation}
\end{lem}
%The simplification is the following. In \cite{M} instead of
%$\eta_{n,-n}^{m_{-n,n}-m'_{-n,n}}$ a more difficult expression is
%placed. But it's action onto the maximal vector coincides with the
%action of $\eta_{n,-n}^{m_{-n,n}-m'_{-n,n}}$.

\proof

  To pass from \eqref{stvm} to  \eqref{stvm1}  we first of all change  $Z_{n,-n}(k)$ to  $\eta_{n,-n}$. Indeed, from Lemma \ref{znnac}  it immediately follows that the results of the actions of  $\eta_{n,-n}$ and $Z_{n,-n}(u)$ on the space  $\mathfrak{sp}_{2n-2}$-highest vectors are proportional (the coefficient of proportionality depends on the vector)

Then using the commutation relations we change the order of factors and come to \eqref{stvm1}.

\endproof

Using Lemma \ref{soo} we obtain the following Lemma

\begin{lem}
\label{dlsl}
$\mathfrak{sp}_{2n-2}$-highest vector in the Molev's approach can be given by the formula

\begin{align}
\begin{split}
\label{stvme}
v=const\prod_{i=1}^{n-1} \varepsilon_{-n,-i}^{m'_{-i,n}-m_{-i,n-1}}
\prod_{i=1}^{n} \varepsilon_{n,-i}^{m_{-i,n}-m'_{-i,n}}v_0
\end{split}
\end{align}

\end{lem}

From lemmas above we obtain the following result.

\begin{cor}
\label{slea} In the restriction problems $\mathfrak{gl}_{n+1}\downarrow\mathfrak{gl}_{n-1}$ and $\mathfrak{sp}_{2n}\downarrow\mathfrak{sp}_{2n-2}$ in the Zhelobenko's approach
the corresponding highest vector is a polynomial in  $a_{-n,...,-k}$, $a_{-n,...-k-1,-1}$, $a_{-n,...,-k-1,1}$, $a_{-n,...-k-2,-1,1}$.

And in the Molev's approach - a polynomial in $a_{-n,...,-k}$, $a_{-n,...-k-1,n}$, $a_{-n+1,...,-k+1}$, $a_{-n+1,...-k,n}$.
\end{cor}

\section{Semimaximal diagrams}

\begin{defn}
We that that the tableau is {\it semimaximal   in the Molev's
approach} if it is of type

\begin{align}
\begin{split}
\label{pmdm} & m_{-n,n}\,\,\,\,\,
\,\,\,\,\,m_{-n+1,n}\,\,\,\,\,\,\,\,\,\,
m_{-n+2,n}\,\,\,\,\,\,\,\,\,\,...\,\,\,\,\,
\,\,\,\,\,m_{-1,n}\,\,\,\,\,\,\,\,\,\,0\\
&\,\,\,\,\,\,\,\,\,\,\,\,\,\,\,\,\,\,\,\,
m'_{-n,n}\,\,\,\,\,\,\,\,\,\,m'_{-n+1,n}\,\,\,\,\,\,\,\,\,\,m'_{-n+2,n}\,\,\,\,\,\,\,\,\,\,...\,\,\,\,\,\,\,\,\,\,m'_{-1,n}\\
&\,\,\,\,\,\,\,\,\,\,\,\,\,\,\,\,\,\,\,\,
\,\,\,\,\,\,\,\,\,\,\,\,\,\,\,\,\,\,\,\,\,\,\,\,\,\,\,\,m'_{-n+1,n}\,\,\,\,\,\,\,\,\,\,m'_{-n+2,n}\,\,\,\,\,\,\,\,\,\,...\,\,\,\,\,\,\,\,\,\,m'_{-1,n},\\
\end{split}
\end{align}
that is if $m'_{-i+1,n}=m_{-i,n-1}$.

\end{defn}

\begin{defn}
We that that the tableau is {\it semimaximal   in the Zhelobenko's
approach} if it is of type

\begin{align}
\begin{split}
\label{pmdzh} & m_{-n,n}\,\,\,\,\,
\,\,\,\,\,m_{-n+1,n}\,\,\,\,\,\,\,\,\,\,
m_{-n+2,n}\,\,\,\,\,\,\,\,\,\,...\,\,\,\,\,
\,\,\,\,\,m_{-1,n}\,\,\,\,\,\,\,\,\,\,0\\
&\,\,\,\,\,\,\,\,\,\,\,\,\,\,\,\,\,\,\,\,
m'_{-n,n}\,\,\,\,\,\,\,\,\,\,m'_{-n+1,n}\,\,\,\,\,\,\,\,\,\,m'_{-n+2,n}\,\,\,\,\,\,\,\,\,\,...\,\,\,\,\,\,\,\,\,\,m'_{-1,n}\\
&\,\,\,\,\,\,\,\,\,\,\,\,\,\,\,\,\,\,\,\,
m'_{-n,n}\,\,\,\,\,\,\,\,\,\,m'_{-n+1,n}\,\,\,\,\,\,\,\,\,\,m'_{-n+2,n}\,\,\,\,\,\,\,\,\,\,...\,\,\,\,\,\,\,\,\,\,\\
\end{split}
\end{align}
that is if $m'_{-i,n}=m_{-i,n-1}$.

\end{defn}

Let us be give a tableau of type \eqref{dia}. Then the corresponding
$\mathfrak{sp}_{2n-2}$-highest vector in the Molev's approach is of
type

\begin{equation}
\label{vm} (m)^M:=const  \prod_{i=1}^{n-1}
\eta_{n,i}^{m_{-i,n-1}-m'_{-i+1,n}} v_1,
\end{equation}

where $v_1$ - is a vector of the semimaximal tableau that is
obtained by the change of the lower row  $m_{-i,n-1}\mapsto
m'_{-i+1,n}$.

From Lemma \ref{soo} it follows that

\begin{equation}
\label{vm1} (m)^M:=const  \prod_{i=1}^{n-1}
\varepsilon_{n,i}^{m_{-i,n-1}-m'_{-i+1,n}} v_1,
\end{equation}

Indeed as it is shown in the next section the action of
$\varepsilon_{n,i}$ and $\eta_{n,i}$  differs only by a constant.

An analogous formula is true in the Zhelobenko's approach, the
$\mathfrak{sp}_{2n-2}$-highest vector can be represented as follows
\begin{equation}
\label{vzh1} (m)^{Zh}:=const  \prod_{i=1}^{n-1}
e_{-1,-i}^{m_{-i,n-1}-m'_{-i+1,n}} v_2,
\end{equation}

where $v_2$ - is a vector of the semimaximal tableau that is
obtained by the change of the lower row  $m_{-i,n-1}\mapsto
m'_{-i,n}$.

From Lemma \ref{soo} it follows that

\begin{equation}
\label{vzh} (m)^{Zh}:=const  \prod_{i=1}^{n-1}
\zeta_{-1,-i}^{m'_{-i,n}-m_{-i,n-1}} v_2,
\end{equation}

\section{The main theorem}

\begin{defn}
Let us be given a tableau  $(m)$ of type \eqref{dia}, define the tableau $(\bar{m})$ by formula

\begin{align}
\begin{split}
\label{diamm}
&m_{-n,n}\,\,\,\,\,\,\,
m_{-n+1,n}\,\,\,\,\,\,\,\,\,\,\,\,\,\,...\,\,\,\,\,\,\,m_{-1,n}\,\,\,\,\,\,\, 0\\
&\,\,\,\,\,\,\,\,\,\,m'_{-n,n}\,\,\,\,\,\,\,m'_{-n+1,n}\,\,\,\,\,\,\,...\,\,\,\,\,\,\,\,\,\,\,m'_{-1,n}\\
&\,\,\,\,\,\,\,\,\,\,\,\,\,\,\,\,\,\,\,\,\,\,\,\,\bar{m}_{-n,n-1}\,\,\,\,\,\,...\,\,\,\,\,\,\,\,\,\,\,\,\,\,\bar{m}_{-2,n-1}\\
\end{split}
\end{align}

 Where $\bar{m}_{-i,n-1}=m'_{-i,n}-(m_{-i,n-1}-m'_{-i+1,n})$.
 \end{defn} Note that if a tableau  $(m)$  is semimaximal in the Zhelobenko's sense that  $(\bar{m})$ is semimaximal in the Molev's sense.  Also  $(\bar{\bar{m}})=(m)$

%Let  $(m)$be a tableau of type  \eqref{dia},  that denotes a
%$\mathfrak{sp}_{2n-2}$-highest vector. The corresponding
%$\mathfrak{sp}_{2n-2}$-highest vector in the Zhelobenko approach we
%denote as   $(m)^{Zh}$,  and the corresponding vector in the Molev's
%approach as $(m)^{M}$.

\begin{thm}
\label{osnvteor}
Under the correspondence $(m)^{Zh}\leftrightarrow (\bar{m})^{M}$ the actions of
$Z^{Zh}_{s}(\mathfrak{sp}_{2n},\mathfrak{sp}_{2n-2})$ and
$Z^{Zh}_{s}(\mathfrak{sp}_{2n},\mathfrak{sp}_{2n-2})$  are conjugated modulo a constant  by the ruler

\begin{align*}
\begin{split}
& \zeta_{1,-k}\sim \eta_{n,-(n-k+1)},\,\,\,\,\,\, \zeta_{-k,1}\sim \eta_{-(n-k+1),n}\\
& \zeta_{-1,-k}\sim \eta_{n-k+1,n},\,\,\,\,\,\, \zeta_{-k,-1}\sim \eta_{n,n-k+1}\\
\end{split}
\end{align*}
\end{thm}

%\begin{cor}
%Совпадают редуцированные матричные элементы  \begin{align*}
%\begin{split}
%& \langle (\bar{m})^{Zh}|F_{1,-k}|(m)^{Zh}\rangle=\langle (\bar{m})^{M}|F_{n,-k}|(m)^{M}\rangle %,\,\,\,\langle (\bar{m})^{Zh}|F_{-k,1}|(m)^{Zh}\rangle=\langle (\bar{m})^{M}|F_{-k,n}|%(m)^{M}\rangle \\
%& \langle (\bar{m})^{Zh}|F_{-1,-k}|(m)^{Zh}\rangle=\langle (\bar{m})^{M}|F_{n,k}|(m)^{M}\rangle %,\,\,\,\langle (\bar{m})^{Zh}|F_{-k,-1}|(m)^{Zh}\rangle=\langle (\bar{m})^{M}|F_{k,n}|%(m)^{M}\rangle.
%\end{split}
%\end{align*}
%\end{cor}

\subsection{The proof in the case of semimaximal tableaux}
Let us prove a statement analogous to the main statement of Section \ref{sviaz}.

\begin{lem}
\label{lm}
Let us be given the algebra $\mathfrak{gl}_{n+1}$ acting in the space with coordinates
$-n,-n+1,...,-1,n$, and  let us be given a subalgebra$\mathfrak{gl}_{n-1}$ that preserves coordinates
$-n+1,...,-1$. Also let $\mathfrak{sp}_{2n}$ act in the space with
coordinates $-n,-n+1,...,n-1,n$, and let $\mathfrak{sp}_{2n-2}$
preserve the coordinates $-n+1,...,n-1$.

Then the problems of
restriction $\mathfrak{gl}_{n+1}\downarrow\mathfrak{gl}_{n-1}$ and
$\mathfrak{sp}_{2n}\downarrow\mathfrak{sp}_{2n-2}$ are equivalent.
\end{lem}
\proof

Indeed take an embedding of $\mathfrak{gl}_{n+1}$ into
$\mathfrak{gl}_{2n}$ as a subalgebra that preserve coordinates
$-n,-n+1,...,-1,n$. Then $\mathfrak{gl}_{n-1}$-highest vectors are
presented as polynomials in variables  $z_{-n,-n+1},...,z_{-n,-1}$,
and $z_{-n,n},z_{-n+1,n},...,z_{-1,n}$.

But in the problem of restriction
$\mathfrak{sp}_{2n}\downarrow\mathfrak{sp}_{2n-2}$ the
$\mathfrak{sp}_{2n-2}$-highest vectors are presented as polynomials
in the same variables.

Also the indicator systems for these polynomials  in the problems of
restriction $\mathfrak{gl}_{n+1}\downarrow\mathfrak{gl}_{n-1}$ and
$\mathfrak{sp}_{2n}\downarrow\mathfrak{sp}_{2n-2}$ coincide.

\endproof

 Using  this Lemma and Lemma \ref{dlsl}
we come to an important conclusion that the Molev's base can be obtained by a procedure analogous to the procedure that defines the Zhelobenko's base but using another restriction  $\mathfrak{sp}_{2n}\downarrow\mathfrak{sp}_{2n-2}$.

\begin{cor}
\label{mcor}
Under the equivalence of restriction problems established in Corollary \ref{sl} the Gelfand-Tsetlin base for $\mathfrak{gl}_{n+1}\downarrow\mathfrak{gl}_{n-1}$  is mapped modulo a constant  to the Molev's base in the restriction problem $\mathfrak{sp}_{2n}\downarrow\mathfrak{sp}_{2n-2}$.

\end{cor}

\begin{lem}
\label{l8}
In the realization of the representation in the function on the
group $Sp_{2n}$ in the Zhelobenko's approach a base vector
corresponding to \eqref{pmdzh} is a polynomial
 $$const \prod_{k=1}^{n-1}   a_{-n,...,-k+1,1}^{m_{-(n-k),n}-m'_{-(n-k),n}}a_{-n,...,-k}^{m'_{-(n-k),n}-m_{-(n-k)+1,n}},$$
and in the Molev's approach the vector corresponding to\eqref{pmdm}
is the polynomial
  $$const \prod_{k=1}^{n-1}  a_{-n,...,-k+1,n}^{m_{-k,n}-m'_{-k,n}}a_{-n,...,-k}^{m'_{-k,n}-m_{-k,n}}.$$
\end{lem}

\proof
The statement of the Lemma in the case of Zhelobenko's
approach  follows from the formula of $\mathfrak{gl}_{n}$-highest
vectors in a $\mathfrak{gl}_{n+1}$-representation  (see \cite{zh2}).
.

Let us prove the Lemma in the case of Molevs's approach.
%The same
%arguments as in Section \ref{sviaz} prove that the problems of
%restriction $\mathfrak{gl}_{n+1}\downarrow\mathfrak{gl}_{n-1}$ and
%$\mathfrak{sp}_{2n}\downarrow\mathfrak{sp}_{2n-2}$ are equivalent.
%Here $\mathfrak{gl}_{n+1}$ acts in the space with coordinates
%$-n,-n+1,...,-1,n$, and $\mathfrak{gl}_{n-1}$ preserves coordinates
%$-n+1,...,-1$. Also $\mathfrak{sp}_{2n}$ acts in the space with
%coordinates $-n,-n+1,...,n-1,n$, and $\mathfrak{sp}_{2n-2}$
%preserves coordinates $-n+1,...,n-1$.
Using
Lemma \ref{soo} we obtain that the vector corresponding to a
semimaxiaml diagram is given by formula

$$
const\prod_{i=1}^{n-1}\varepsilon_{n,-i}^{m_{i,n}-m'_{i,n}}v_0.
$$

% In both cases
%$$v=\prod_{i=n}^1a_{-n,...,-i}^{m_{-i,n}-m_{-i+1,n}}$$ is a maximal vector.

Using the formulas for the action of
$\varepsilon_{n,i}$ onto the determinants we immediately obtain the
statement of the Lemma.

\endproof

\begin{cor}
\label{sl}
In the realization of the representation in the function on the
group $Z$ in the Zhelobenko's
 approach a base vector corresponding to \eqref{pmdzh} is a polynomial  $$const \prod_{k=1}^{n-1} z_{-k,1}^{m_{-(n-k),n}-m'_{-(n-k),n}},$$
and in the Molev's approach the vector corresponding to\eqref{pmdm}
is the polynomial  $$const \prod_{k=1}^{n-1}
z_{-(n-k+1),n}^{m_{-k,n}-m'_{-k,n}}.$$
\end{cor}

%Let us reformulate this lemma on the language of function on the whole group.

%This corollary gives a formula for  $v_1$  in the formulas \eqref{vm}, \eqref{vm2}.

Now let us return to the proof of the Theorem in the case of semi-maximal tableaux.
For semimaximal tableaux the correspondence  $(m)^{Zh}\leftrightarrow (\bar{m})^M$ means the correspondence

\begin{align}
\begin{split}
\label{alr} &a_{-n,...,-k}\leftrightarrow  a_{-n,...,-k}\\
&a_{-n,...,-k+1,1}\leftrightarrow  a_{-n,...,-n+k,n}
\end{split}
\end{align}

Thus we need to prove the following statement.

\begin{lem}
 The actions of  $\zeta_{1,-k}$ and $\eta_{n,-(n-k+1)}$, $\zeta_{-k,1}$
 and
$\eta_{-(n-k+1),n}$, $k=n,...,1$  on the semimaximal tableaux under
the correspondence \eqref{alr} are conjugated modulo multiplication
by a constant.
\end{lem}

\proof

It is sufficient to prove that the actions of
 $\zeta_{1,-k}$ and $\eta_{n,-(n-k)}$, $\zeta_{-k,1}$
 and
$\eta_{-(n-k+1),n}$  on $a_{-n,...,-p}$, $a_{-n,...,-p-1,1}$ and
$a_{-n,...,-p}$, $a_{-n,...,-p-1,n}$ are conjugated under the
correspondence \eqref{alr}

Under the action of $\zeta_{1,-k}$ the determinant are changed as
follows

\begin{align*}
&a_{-n,...,-p}\mapsto \begin{cases} 0,\,\,\, p\neq k\\ a_{-n,...,-p-1,1},\,\,\, p= k \end{cases}\\
&a_{-n,...,-p-1,1}\mapsto 0
\end{align*}

 Under the action of  $\zeta_{-k,1}$   the determinant are changed as follows

\begin{align*}
&a_{-n,...,-p}\mapsto  0 \\
&a_{-n,...,-p-1,1}\mapsto \begin{cases} 0,\,\,\, p\neq k\\ a_{-n,...,-p-1,-p},\,\,\, p=k \end{cases}
\end{align*}

 Under the acion of  $\eta_{n,-(n-k+1)}$   the determinant are changed as follows

\begin{align*}
&a_{-n,...,-p}\mapsto \begin{cases} 0,\,\,\, p\neq n-k+1\\ a_{-n,...,-p-1,n},\,\,\, p\neq n-k+1 \end{cases}\\
&a_{-n,...,-p-1,n}\mapsto 0
\end{align*}

Under the acion of  $\eta_{-(n-k+1),n}$   the determinant are changed as follows

\begin{align*}
&a_{-n,...,-p}\mapsto  0\\
&a_{-n,...,-p-1,n}\mapsto \begin{cases} 0,\,\,\, p\neq n-k+1\\ a_{-n,...,-p-1,-p},\,\,\, p\neq n-k+1 \end{cases}
\end{align*}

If one compares the corresponding formulas, then using Lemma \ref{l8} one  obtains the
conclusion of the Lemma.

\endproof

\subsection{Thee proof in the case of an arbitrary tableau}

We use the formulas  \eqref{vm} and \eqref{vzh}.

It is sufficient to check that the following actions are conjugated

\begin{align*}
& \eta_{n-k+1,n}(\bar{m})^M\leftrightarrow \zeta_{-1,-k}(m)^{Zh},\,\,\, \eta_{n,n-k+1}(\bar{m})^M\leftrightarrow \zeta_{-k,-1}(m)^{Zh} \\
& \eta_{n,-(n-k+1)}(\bar{m})^M\leftrightarrow \zeta_{-k,n}(m)^{Zh},\,\,\,\eta_{-(n-k+1),n}(\bar{m})^M\leftrightarrow \zeta_{n,-k}(m)^{Zh}.\\
\end{align*}

\subsubsection{The proof for  $\zeta_{-1,-k} $ and  $\eta_{n-k+1,n}$}  Since the operators $\zeta_{-1,-i}$ commute,
the multiplication of  \eqref{vzh} onto  $\zeta_{-1,-k}$ means that
we make the following change   $$m'_{-k,n}-m_{-k,n-1}\mapsto
m'_{-k,n}-m_{-k,n-1}+1.$$  Analogously the multiplication of
\eqref{vm} onto $\eta_{n-k+1,n}$ is equivalent to the change
$$m'_{-k,n}-m_{-k,n-1}\mapsto m'_{-k,n}-m_{-k,n-1}+1.$$ Thus both
operations make a change

\begin{equation*}
m_{-k,n-1}\mapsto m_{-k,n-1}-1,
\end{equation*}

Hence these operations are conjugate

\subsubsection{The proof for $\zeta_{-k,-1} $ and  $\eta_{n,n-k+1}$}

In this case we must compare the results of multiplication of
\eqref{vzh} on  $\zeta_{-k,-1}$ and the result of multiplication of
\eqref{vm} on $\eta_{k,n}$.

 The operators  $\zeta_{-1,-i}$,  $i>k$ and  $\zeta_{-k,-1}$commute and $\zeta_{-k,-1}\zeta_{-1,-k}=1$. Thus the multiplication of
 \eqref{vzh} onto  $\zeta_{-k,-1}$ means that we make a change  $$m'_{-k,n}-m_{-k,n-1}\mapsto m'_{-k,n}-m_{-k,n-1}-1.$$
  analogously the multiplication of \eqref{vm} onto  $\eta_{n,n-k+1}$ is equivalent to the change
    $$m'_{-k,n}-m_{-k,n-1}\mapsto m'_{-k,n}-m_{-k,n-1}-1.$$ Thus both operations make a change

\begin{equation*}
m_{-k,n-1}\mapsto m_{-k,n-1}+1,
\end{equation*}

Hence these operations are conjugate.

\subsubsection{The prof for $\zeta_{1,-k} $ and  $\eta_{-(n-k+1),n}$.}

We must compare the results of multiplication of \eqref{vzh} onto
$\zeta_{1,-k}$ and the result of multiplication \eqref{vm} onto
$\eta_{-(n-k+1),n}$.

Consider the first expression.   Considering the defining commutation relations for the algabra
  $Z^{zh}_s(\mathfrak{sp}_{2n},\mathfrak{sp}_{2n-2})$ one obtains a relation   $$\zeta_{1,-k}\zeta_{-1,-k+1}=\zeta_{-1,-k}\zeta_{1,-k+1}.$$

Let in the formula \eqref{vzh} the factor $\zeta_{-1,-k+1}$  enters
with a nonzero exponent  $b=m'_{-(k+1),n}-m_{-(k+1),n-1}> 0 $, then
\begin{align}
\begin{split}
\label{fo1}
&\zeta_{1,-k}(m)^{Zh}=\zeta_{1,-k}(...\zeta_{-1,-k}^{a}\zeta_{-1,-k+1}^{b}...)v_1=(...\zeta_{-1,-k}^{a}\zeta_{1,-k}\zeta_{1,-k+1}^{b}...)v_1=\\
&=(...\zeta_{-1,-k}^{a}\zeta_{1,-k}\zeta_{-1,-k}\zeta_{-1,-k}^{b-1}...)v_1
=(...\zeta_{-1,-k}^{a}\zeta_{-1,-k}\zeta_{1,-k+1}\zeta_{-1,-k+1}^{b-1}...)v_1=\\
&=(...\zeta_{-1,-k}^{a+1}\zeta_{-1,-k+1}^{b-1}...)\zeta_{-1,-k+1}v_1
\end{split}
\end{align}

And if   $b=0 $, then

\begin{align}
\begin{split}
\label{fo2}
&\zeta_{1,-k}(m)^{Zh}=\zeta_{1,-k}(...)v_1=(...)\zeta_{1,-k}v_1
\end{split}
\end{align}

Also one has in  $Z^{M}_s(\mathfrak{sp}_{2n},\mathfrak{sp}_{2n-2})$
a relation   $$\eta_{-k,n}\eta_{n,k+1}=\eta_{n,k}\eta_{-k-1,n}.$$ As a corollary we get a relation
$$\eta_{-(n-k+1),n}\eta_{n,(n-k+2)}=\eta_{n,(n-k+1)}\eta_{-(n-k+2),n}.$$

Let in the formula \eqref{vm} the factor $\eta_{n,k}$ enters with a
nonzero exponent $b=m'_{-(k-1),n}-m_{-(k-1),n-1}> 0 $, then
\begin{align}
\begin{split}
\label{for1}
&\eta_{-(n-k+1),n}(\bar{m})^{M}=\eta_{-(n-k+1),n}(...\eta_{n,(n-k+1)}^{a}\eta_{n,(n-k+2)}^{b}...)v_2=\\&=(...\eta_{n,(n-k+1)}^{a}\eta_{-(n-k+1),n}\eta_{n,(n-k+2)}^{b}...)v_2=\\
&=(...\eta_{n,(n-k+1)}^{a}\eta_{-(n-k+1),n}\eta_{n,(n-k+2)}\eta_{n,(n-k+2)}^{b-1}...)v_2=\\&=
=(...\eta_{n,(n-k+1)}^{a}\eta_{n,(n-k+1)}\eta_{-(n-k+2),n}\eta_{n,(n-k+2)}^{b-1}...)v_2=\\
&=(...\eta_{n,(n-k+1)}^{a+1}\eta_{n,(n-k+2)}^{b-1}...)\eta_{-(n-k+2),n}v_2
\end{split}
\end{align}

And if   $b$, then

\begin{align}
\begin{split}
\label{for2}
&\eta_{-(n-k+1),n}(\bar{m})^{M}=\eta_{-(n-k+1),n}(...)v_2=(...)\eta_{-(n-k+1),n}v_2
\end{split}
\end{align}

Compare the formulas  \eqref{fo1} and  \eqref{for1},  \eqref{fo2} and  \eqref{for2}.
  We come to the following result. If  $m'_{-(k-1),n}-m_{-(k-1),n-1}\neq 0$,
  then in both cases the following changes occur. Firstly $m_{-(k-1),n-1}\mapsto m_{-(k-1),n-1}+1$,
  secondly onto  $v_1$ and $v_2$  there acts  $\zeta_{1,-k+1}$  or  $\eta_{-(n-k+2),n}$
   correspondingly.   And if  $m'_{-(k-1),n}-m_{-(k-1),n-1}=0$, then onto  $v_1$ and $v_2$ there acts  $\zeta_{1,-k}$  or  $\eta_{-(n-k+1),n}$ correspondingly.

Hence both actions are conjugated modulo multiplication by a
constant under the considered correspondence.

\subsubsection{The proof for $\zeta_{-k,1} $ and  $\eta_{n,-(n-k+1)}$.}

We must compare the results of multiplication of \eqref{vzh} onto
$\zeta_{-k,1} $ and the result of multiplication of \eqref{vm} onto
$\eta_{n,-k}$.

Consider the first expression, one in  has the simplified algebra of Mikelsson-Zhelobenko
 $Z^{Zh}_s(\mathfrak{sp}_{2n},\mathfrak{sp}_{2n-2})$   a relation $$\zeta_{-k,1}\zeta_{-1,-k}=\zeta_{-1,-(k-1)}\zeta_{-(k-1),1}.$$
Let in the formula \eqref{vzh} the factor $\zeta_{1,k+1}$  occurs
with a nonzero exponent  $b=m'_{-k,n}-m_{-k,n-1}> 0 $ then
\begin{align}
\begin{split}
\label{foo1}
&\zeta_{-k,1}(m)^{Zh}=\zeta_{-k,1}(...\zeta_{-1,-k}^{b}\zeta_{-1,-(k-1)}^{a}...)v_1=(...\zeta_{-k,1}\zeta_{-1,-k}^{b}\zeta_{-1,-(k-1)}^{a}...)v_1=\\
&=(...\zeta_{-k,1}\zeta_{-1,-k}\zeta_{-1,-k}^{b-1}\zeta_{-1,-(k-1)}^{a}...)v_1
=(...\zeta_{-1,-(k-1)}\zeta_{-(k-1),1}\zeta_{-1,-k}^{b-1}\zeta_{-1,-(k-1)}^{a}...)v_1=\\
&=(...\zeta_{-1,-k}^{b-1}\zeta_{-1,-(k-1)}^{a+1}...)\zeta_{-(k-1),1}v_1
\end{split}
\end{align}

And if  $b=m'_{-k,n}-m_{-k,n-1}=0 $ then

\begin{align}
\begin{split}
\label{foo2}
&\zeta_{-k,1}(m)^{Zh}=\zeta_{-k,1}(...)v_1=(...)\zeta_{-k,1}v_1
\end{split}
\end{align}

From the other hand one has in $Z^{M}_s(\mathfrak{sp}_{2n},\mathfrak{sp}_{2n-2})$   a relation  $$\eta_{n,-k}\eta_{n,k}=\eta_{n,k+1}\eta_{n,-(k+1)}.$$ As a corollary we obtain a relation
 $$\eta_{n,-(n-k+1),}\eta_{n,(n-k+1)}=\eta_{n,(n-k+2)}\eta_{n,-(n-k+2)}.$$

%$f_{n,k+1}f_{n,-k-1}=f_{n,-k}f_{n,k}$

Let in the formulas \eqref{vm} the factor  $\eta_{n,(n-k+1)}$  occurs with a nonzero exponent  $b=m'_{-k,n}-m_{-k,n-1}> 0 $ then
\begin{align}
\begin{split}
\label{foor1}
&\eta_{n,-(n-k+1)}(\bar{m})^{M}=\eta_{n,-(n-k+1)}(...\eta_{n,(n-k+1)}^{b}\eta_{n,(n-k+2)}^{a}...)v_2=\\&=(...\eta_{n,-(n-k+1)}\eta_{n,(n-k+1)}^{b}\eta_{n,(n-k+2)}^{a}...)v_2=\\
&=(...\eta_{n,-(n-k+1)}\eta_{n,(n-k+1)}\eta_{n,(n-k+1)}^{b-1}\eta_{n,(n-k+2)}^{a}...)v_2
=\\&=(...\eta_{n,(n-k+2)}\eta_{n,-(n-k+2)}\eta_{n,(n-k+1)}^{b-1}\eta_{n,(n-k+2)}^{a}...)v_2=\\
&=(...\eta_{n,(n-k+1)}^{b-1}\eta_{n,(n-k+2)}^{a+1}...)\eta_{n,-k+1}\eta_{n,(n-k+2)}v_2
\end{split}
\end{align}

If  $b=m'_{-k,n}-m_{-k,n-1}=0 $ then

\begin{align}
\begin{split}
\label{foor2}
&\eta_{n,-(n-k+1)}(\bar{m})^{M}=\eta_{n,-(n-k+1)}(...)v_2=(...)\eta_{n,-(n-k+1)}v_2
\end{split}
\end{align}

Compare the formulas  \eqref{foo1} and  \eqref{foor1},  \eqref{foo2} and  \eqref{foor2}.  We come to the following result. If  $m'_{-k,n}-m_{-k,n-1}\neq 0$, then in both cases the following changes occur. Firstly  $m_{-k,n-1}\mapsto m_{-k,n-1}+1$,  secondly onto  $v_1$ and $v_2$  there acts $\zeta_{1,-(k-1)}$  or   $\eta_{n,-(n-k+2)}$  correspondingly.   And if   $m'_{-k,n}-m_{-k,n-1}= 0$, , then onto  $v_1$ and $v_2$ there acts  $\zeta_{1,-k}$  or   $\eta_{n,-(n-k+1)}$ correspondingly.

Hence both actions are conjugated modulo multiplication by a constant under the considered correspondence.

\end{document}